\documentclass{article}
\usepackage{amsmath}
\usepackage{graphicx}
\usepackage{latexsym}
\usepackage{hyperref}

\usepackage{amsfonts}
\usepackage[all]{xy}
\newcommand{\R}{\mathbb{R}}

\def\p{\partial}
\def\n{\nabla}

\def\a{\alpha}
\def\b{\beta}

\newcommand{\ri}{\rightarrow}
\newcommand{\na}{\nabla}
\newcommand{\pa}{\partial}
\newcommand{\no}{\nonumber}

\newcommand{\ra}{\rangle}
\newcommand{\fr}{\frac}

\newcommand\be{\begin{eqnarray}}
\newcommand\en{\end{eqnarray}}
\def\en{\end{eqnarray}}

\def\H{\mathbf{H}}

\def\div{{\mathrm{div}}}

\def\L{\overline{\Delta}}

\def\ino{\int_\Omega}
\def\inpo{\int_{\partial\Omega}}
\def\D{\Delta}

\def\la{\lambda}

\def\om{\Omega}

\def\({\left(}
\newcommand{\Om}{\Omega}
\newcommand{\laa}{\lambda}
\newcommand{\ov}{\overline}

\def\){\right)}

\newtheorem{theorem}{Theorem}[section]
\newtheorem{lemma}[theorem]{Lemma}

\newtheorem{corollary}[theorem]{Corollary}

\newtheorem{remark}[theorem]{Remark}

\begin{document}
\setcounter{page}{1}
\title{Isoperimetric bounds for eigenvalues of the Wentzell-Laplace, the Laplacian and a biharmonic Steklov problem}
\author{Feng Du$^{a}$, Jing Mao$^{b}$,
Qiaoling Wang$^{c}$, Changyu Xia$^{c}$}
\date{}
\protect\footnotetext{\!\!\!\!\!\!\!\!\!\!\!\!{ MSC 2020: 35P15;
53C40; 58C40.}
\\
{ ~~Key Words: The Wentzell-Laplace operator, biharmonic Steklov
problem, eigenvalues, Hadamard manifolds.} }
\maketitle ~~~\\[-15mm]

\begin{center}
{\footnotesize  $a$.  School of Mathematics and Physics Science,\\
Jingchu University of Technology, \\
Jingmen, 448000, China\\
$b$. Faculty of Mathematics and Statistics,\\
Key Laboratory of Applied Mathematics of Hubei Province,\\
 Hubei University, Wuhan, 430062, China \\
 $c$. Departamento de Matem\'atica, Universidade de Brasilia,\\
70910-900-Brasilia-DF, Brazil\\
 Emails:  defengdu123@163.com (F. Du), jiner120@163.com (J. Mao),\\
wang@mat.unb.br (Q. Wang), xia@mat.unb.br (C. Xia).  }
\end{center}


\begin{abstract}
In this paper, we prove some isoperimetric bounds for lower order
eigenvalues of the Wentzell-Laplace operator on bounded domains of a
Euclidean space or a Hadamard manifold, of the Laplacian on closed
hypersurfaces of a Euclidean space or a Hadamard manifold, and of a
biharmonic Steklov problem on bounded domains of a Euclidean space.
Especially, interesting rigidity results can be obtained if sharp
bounds were achieved.
 \end{abstract}

\markright{\sl\hfill  F. Du, J. Mao, Q. Wang, C. Xia \hfill}

\section{Introduction}
\renewcommand{\thesection}{\arabic{section}}
\renewcommand{\theequation}{\thesection.\arabic{equation}}
\setcounter{equation}{0} \label{intro}

Let $M$ be an $n$-dimensional compact Riemannian manifold with
boundary\footnote{For the eigenvalue problem (\ref{int1}), the
regularity assumption for $\partial M$ should be made such that the
embedding $H(M)\subset L^{2}(M)\cup L^{2}(\partial M)$ is compact,
which is the essential requirement so that the Laplacian in
(\ref{int1}) has a discrete spectrum. Here, $H(M)$ defined by
(\ref{FC}) below is the admissible space of the eigenvalue problem
(\ref{int1}). In fact, for the eigenvalue problem (\ref{int1}),
 the regularity assumption ``\emph{$\partial M$ is Lipschitz continuous}" can already make sure that the embedding $H(M)\subset L^{2}(M)\cup L^{2}(\partial M)$
 is compact. Generally speaking, different eigenvalue problems might have different regularity requirements on the boundary (if nonempty) to have discrete spectrum.
 In order to avoid the repeated argument on regularity of the boundary (which is not the main part of this paper), without specification we always assume that the
 boundary regularity is \emph{smooth} for all eigenvalue
 problems considered in this paper. Moreover, this smooth
 assumption would be mentioned in the sequel if necessary.
 } $\partial M$
and denote by $\overline{\Delta}$ and $\Delta$ the Laplace-Beltrami
operators on $M$ and $\partial M$, respectively. Assume that $\beta$
is a real number and consider the following eigenvalue problem with
the Wentzell boundary condition
\be\label{int1}\left\{\begin{array}{l} \ov{\Delta} u =0 \ \ \ \ \ \
\ \ \ \ \ \ \ \ \ \ \ \ \ \ \ \ {\rm in \ }\ M,\\ -\beta\Delta
u+\pa_{\nu} u= \laa u\ \ \ \quad~ {\rm on  \ } \pa M,
\end{array}\right.
\en where  $\pa_{\nu}$ denotes the derivative w.r.t. the outward
unit normal vector $\nu$ to $\partial M$. This problem has already
been studied in \cite{DKL,G} when $M$ is a bounded domain in a
Euclidean space. Recently, some interesting estimates for the
nonzero eigenvalues of the problem (\ref{int1}) and its weighted
version have been obtained (see, e.g., \cite{DWX,WX3,YWMD}).
Besides, following the convention in \cite{YWMD}, we call
(\ref{int1}) \emph{the Wentzell eigenvalue problem of the
Laplacian}.
 Note
that when $\beta=0$, the eigenvalue problem (\ref{int1}) becomes the
classical Steklov problem
 \be \label{int2}\left\{\begin{array}{l} \ov{\Delta} u =0 \ \
\ \ \ \ ~ {\rm in \ \ } M,\\ \pa_{\nu} u=  p u \ \  ~~ {\rm on \ \ }
\pa M,
\end{array}\right.
\en which has a discrete spectrum\footnote{Of course, each element
in the discrete spectrum is called \emph{eigenvalue}.} consisting in
a non-decreasing sequence\footnote{Clearly, eigenvalues $p_{i}$ in
the sequence (\ref{SES}) should be written as $p_{i}(M)$ for
accuracy. However, for convenience, in the sequel, we prefer to
simplify the notations for eigenvalues (of different type) discussed
in this paper, that is, we separately write $p_{i}(M)$,
$\lambda_{i,\beta}(\Omega)$, $\lambda_{i}(M)$,
$\lambda_{i,\tau}(\Omega)$ as $p_{i}$, $\lambda_{i,\beta}$,
$\lambda_{i}$ and $\lambda_{i,\tau}$. We also make an agreement that
these notations would be written completely if necessary.}
 \be \label{SES}
p_0=0<p_1\leq p_2\leq\cdots \rightarrow +\infty.\en For the Steklov
eigenvalue problem (\ref{int2}), many interesting results have been
obtained and one can find some of them -- see, e.g.,
\cite{ceg,E1,E2,FS,HPS,kk, K,ST,WX1,WX2,YWMD,yy} and the references
therein.

When $\beta\geq 0$, the spectrum of the Laplacian with Wentzell
boundary condition, i.e. the spectrum of the eigenvalue problem
(\ref{int1}), consists in an non-decreasing countable sequence of
eigenvalues \be \no\laa_{0, \beta}=0<\laa_{1, \beta}\leq\laa_{2,
\beta}\leq\cdots\rightarrow +\infty, \en with corresponding real
orthonormal (in $L^2(\pa M)$ sense) eigenfunctions $u_0, u_1,
u_2,\cdots.$ We adopt the convention that each eigenvalue is
repeated according to its multiplicity.  Consider the Hilbert space
\be \label{FC} H(M)=\{ u\in H^1(M), {\rm  Tr}_{\pa M}(u)\in H^1(\pa
M) \}, \en where ${\rm Tr}_{\pa M}$ is the trace operator on
$\partial M$. We define on $H(M)$ the two bilinear forms \be
\label{TF} A_{\beta}(u, v)=\int_{M}\ov{\nabla}u\cdot\ov{\nabla}v
+\beta\int_{\pa M} \na u\cdot\na v, \quad \ B(u, v)=\int_{\pa M}uv,
\en where, as before, $\ov{\nabla}$ and $\na$ are the gradient
operators on $M$ and $\pa M$, respectively. Here, volume elements in
the above integrals have been dropped.\footnote{For convenience, in
the sequel we will drop the integral measures for all integrals
except it is necessary.} Since we assume  that $\beta$  is
nonnegative, two bilinear forms defined by (\ref{TF}) are positive
and the variational characterization of the $k$-th eigenvalue is \be
\label{vc}\laa_{k, \beta}=\min\left\{\frac{A_{\beta}(u, u)}{B(u,
u)}\Bigg{|} u\in H(M), u\neq 0,\ \int_{\pa M} u u_i=0,
i=0,\cdots,k-1\right\}. \en Clearly, when $k=1$, the minimum is
taken over the functions orthogonal (in $L^2(\pa M)$ sense) to the
eigenfunctions associated to $\laa_{0, \beta}= 0,$ i.e., nonzero
constant functions.

If $M$ is an $n$-dimensional Euclidean ball of radius $R$, then (cf.
\cite{DKL})
 \begin{eqnarray}  \label{add-1}
\laa_{1, \beta}=\laa_{2, \beta}=\cdots \laa_{n,
\beta}=\fr{(n-1)\beta + R}{R^2}
\end{eqnarray}
and  the corresponding eigenfunctions are the coordinate functions
$x_i,  i=1,\cdots,n.$ For the Steklov problem (\ref{int2}), Brock
 \cite{B} proved the following well-known result:

\begin{theorem} (\cite{B}) \label{th0} Let $\Omega$ be a bounded domain with smooth boundary in the Euclidean
$n$-space $\R^n$ and denote by $p_1(\om),\cdots p_n(\om)$ the first
$n$ nonzero Steklov eigenvalues of $\om$. Then  \be\label{th0.0}
\sum_{i=1}^n\fr 1{p_i(\om)} \geq
n\left(\fr{|\om|}{\omega_{n}}\right)^{\fr 1n}, \no\en
  where $\omega_n$ and $|\om|$
denote the volume of the unit ball in $\R^n$ and of $\om$,
respectively.\footnote{In the sequel, without specification, for a
given geometric subject, $|\cdot|$ would denote its Hausdorff
measure.}
\end{theorem}

The proof of Brock's theorem is a nice application of an inequality
for sums of reciprocals of eigenvalues shown by Hile-Xu in \cite{hx}
and a weighted isoperimetric inequality obtained by
Betta-Brock-Mercaldo-Posteraro in \cite{bb}. Brock's method has been
used by Dambrine-Kateb-Lamboley in \cite{DKL} to obtain an estimate
for eigenvalues of the Wentzell eigenvalue problem of the Laplacian
(\ref{int1}).

In the first part of this paper, we use the variational
characterization (\ref{vc}) to get the following Brock-type
isoperimetric bound.

\begin{theorem} \label{th1} Let $\beta\geq 0$ and $\Omega$ be a bounded domain with smooth boundary $\p \Omega$ in
 $\R^n$.  Let $|\p\Omega|$ be the area of $\p\Omega$ and denote by $\laa_{1, \beta}\leq\laa_{2, \beta}\leq\cdots\leq \laa_{n, \beta}$
 the first $n$ nonzero
 eigenvalues of the following eigenvalue problem with the Wentzell boundary condition
\be\label{th1.1}\left\{\begin{array}{l} \ov{\Delta} u =0 \ \ \ \ \ \
\ \ \ \ \ \ \ \ \ \ \ \ \ \ \ \ \ {\rm in  \ }\ \om,\\ -\beta\Delta
u+\pa_{\nu} u= \laa u\ \ \ \ \quad ~~ {\rm on  \ } \pa \om.
\end{array}\right.
\en
Then we have
\begin{small}
\be\label{th1.2}
\sum_{i=1}^n\frac{|\Omega|}{\lambda_{i,\b}}+\sum_{i=1}^{n-1}\frac{\b|\p\Omega|}{\lambda_{i,\b}}
\geq n\omega_n\left(\fr{|\om|}{\omega_{n}}\right)^{1+\fr
1n}\left[1+\fr{(n+1)(2^{1/n}-1)}{4n}\left(\fr{|\om\Delta
B|}{|B|}\right)^2\right],
 \en
\end{small}
where $B$ is the ball of volume $|\om|$ and with the same center of
mass than $\pa \om$, and $\om\Delta B$ is the symmetric difference
of $\om$ and $B$. Equality holds in (\ref{th1.2}) if and only if
$\om$ is a ball.
\end{theorem}

\begin{remark}{\rm Taking $\beta=0$ in  Theorem \ref{th1}, one gets
Theorem \ref{th0} directly.}
\end{remark}

 Using the monotonicity of eigenvalues $\lambda_{i,\beta}$'s and
Theorem \ref{th1} immediately, we get
\be
&&\qquad \qquad \frac{1}{\lambda_{1,\beta}}\left[n|\Omega|+(n-1)\beta|\partial\Omega|\right]\geq\nonumber\\
&&\sum_{i=1}^n\frac{|\Omega|}{\lambda_{i,\b}}+\sum_{i=1}^{n-1}\frac{\b|\p\Omega|}{\lambda_{i,\b}}
\geq n\omega_n\left(\fr{|\om|}{\omega_{n}}\right)^{1+\fr
1n}\left[1+\fr{(n+1)(2^{1/n}-1)}{4n}\left(\fr{|\om\Delta
B|}{|B|}\right)^2\right], \nonumber
 \en
which directly implies the following eigenvalue estimate.

\begin{corollary} \label{coro-2}
 Under the assumptions in Theorem \ref{th1}, we have
 \begin{eqnarray*}
\lambda_{1,\beta}(\Omega)\leq\left[1+\beta\frac{(n-1)|\partial\Omega|}{n|\Omega|}\right]\left(\frac{w_{n}}{|\Omega|}\right)^{\frac{1}{n}}
\left[1+\fr{(n+1)(2^{1/n}-1)}{4n}\left(\fr{|\om\Delta
B|}{|B|}\right)^2\right]^{-1},
 \end{eqnarray*}
where equality holds if and only if $\om$ is a ball.
\end{corollary}

Corollary \ref{coro-2} and the truth (\ref{add-1}), together with
the classical isoperimetric inequality for bounded Euclidean domains
with fixed volume (see, e.g., \cite[Chapter 1]{RS}), imply the
following estimate:

\begin{corollary} \label{coro-3}
Under the assumptions in Theorem \ref{th1}, we have
 \begin{eqnarray*}
\lambda_{1,\beta}(\Omega)\leq
\lambda_{1,\beta}(B)\cdot\frac{n|\Omega|+(n-1)\beta|\partial\Omega|}{n|\Omega|+(n-1)\beta|\partial
B|}\cdot \left[1+\fr{(n+1)(2^{1/n}-1)}{4n}\left(\fr{|\om\Delta
B|}{|B|}\right)^2\right]^{-1},
 \end{eqnarray*}
 where similarly $|\partial B|$ denotes the area of the sphere $\partial B$ (i.e., the boundary of the ball
 $B$ of volume $|\Omega|$). Equality holds if and only if $\Omega$ is
 a ball, and moreover, in this situation,
 $\lambda_{1,\beta}(\Omega)=
\lambda_{1,\beta}(B)=\fr{(n-1)\beta + R}{R^2}$ with
$R=\left(\frac{|\Omega|}{w_{n}}\right)^{1/n}$.
\end{corollary}

\begin{remark} \label{remark-1}
\rm{ (1) If $\beta=0$, then the conclusion of Corollary \ref{coro-3}
degenerates into
\begin{eqnarray*}
p_{1}(\Omega)\leq p_{1}(B)\cdot
\left[1+\fr{(n+1)(2^{1/n}-1)}{4n}\left(\fr{|\om\Delta
B|}{|B|}\right)^2\right]^{-1},
\end{eqnarray*}
where equality holds if and only if $\Omega$ is
 a ball (i.e., in this case, $\Omega$ is the ball $B$). This result
 indicates:
  \begin{itemize}
\item \emph{Among all bounded Euclidean domains of smooth boundary with
fixed volume, the ball maximizes the first nonzero Steklov
eigenvalue}.\footnote{As shown in the footnote of $1^{\mathrm{st}}$
page, the regularity of the boundary can be weaken to Lipschitz
continuity. However, by (1) of Remark \ref{remark-1}, we do not know
whether the ball $B$ is the only domain such that the functional
$p_{1}(\cdot):\Omega\mapsto p_{1}(\Omega)$ attains its maximum value
or not. But Brock
 \cite{B} has given an affirmative answer already -- the ball is the only possibility to get the maximum value.}
  \end{itemize}
This spectral isoperimetric inequality can be obtained directly by
Theorem \ref{th0} and has already been pointed out by Brock
 \cite{B}.
 \\
 (2) When $\beta>0$, by the classical isoperimetric inequality, we know that under the
constraint $|\Omega|=|B|$, one has
\begin{eqnarray*}
1\leq\frac{n|\Omega|+(n-1)\beta|\partial\Omega|}{n|\Omega|+(n-1)\beta|\partial
B|}<1+\beta\frac{(n-1)|\partial\Omega|}{n|\Omega|},
\end{eqnarray*}
where equality holds if and only if $\Omega$ is the ball $B$. Hence,
Corollary \ref{coro-3} gives a partial answer to
Dambrine-Kateb-Lamboley's conjecture proposed in \cite[page
412]{DKL}. \\
 (3) For the eigenvalue problem (\ref{th1.1}), Dambrine, Kateb and Lamboley \cite[Theorem 1.1]{DKL}
 gave
 the estimate:
\begin{small}
\be \label{add-2}
\left(|\Omega|+\beta\Lambda[\Omega]\right)\sum_{i=1}^n\frac{1}{\lambda_{i,\b}}
\geq n\omega_n\left(\fr{|\om|}{\omega_{n}}\right)^{1+\fr
1n}\left[1+\fr{(n+1)(2^{1/n}-1)}{4n}\left(\fr{|\om\Delta
B|}{|B|}\right)^2\right],
 \en
\end{small}
where $\Lambda[\Omega]$ is spectral radius of the symmetric and
positive semidefinite matrix $p(\Omega)=(p_{ij})_{n\times n}$
defined as
 \begin{eqnarray*}
p_{ij}=\int_{\partial\Omega}\left(\delta_{ij}-\nu_{i}\nu_{j}\right)
\end{eqnarray*}
with, similarly, $\nu$ is the outward unit normal vector to
$\partial\Omega$. They \cite[Lemma 2.4]{DKL} also proved that the
matrix $p(\Omega)=(p_{ij})_{n\times n}$ defined as above should be
symmetric, positive definite, and its spectral radius
$\Lambda[\Omega]$ satisfies
\begin{eqnarray} \label{add-3}
(n-1)|\partial\Omega|\geq\Lambda[\Omega]\geq\frac{n-1}{n}|\partial\Omega|.
\end{eqnarray}
In particular, among sets of given volume, the spectral radius is
minimal for the ball.

On the other hand, by direct calculation, one has
\begin{eqnarray} \label{add-4}
\frac{n-1}{n}\leq\frac{\sum_{i=1}^{n-1}\frac{1}{\lambda_{i,\b}}}{\sum_{i=1}^n\frac{1}{\lambda_{i,\b}}}<1,
\end{eqnarray}
with equality if and only if
$\lambda_{1,\beta}=\lambda_{2,\beta}=\cdots=\lambda_{n,\beta}$.
Hence, combining (\ref{add-2}), (\ref{add-3}) with (\ref{add-4}), it
is easy to know that:
\begin{itemize}

\item If
$\Lambda[\Omega]\in\left[\frac{n-1}{n}|\partial\Omega|,|\partial\Omega|\right)$,
for our estimate (\ref{th1.2}) and Dambrine-Kateb-Lamboley's
estimate (\ref{add-2}), one does not know which one is better;

\item If
$\Lambda[\Omega]\in[|\partial\Omega|,(n-1)|\partial\Omega|]$, then
our estimate (\ref{th1.2}) is sharper than Dambrine-Kateb-Lamboley's
estimate (\ref{add-2}).
\end{itemize}
 (4) By (\ref{add-3}), it is easy to know that the estimate in
 Corollary \ref{coro-2} here is sharper than the one given in
 \cite[Corollary 1.2]{DKL}.
}
\end{remark}

For the Wentzell eigenvalue problem of the Laplacian, we also have:

\begin{theorem}\label{th2} Let $Q^n$ be a  Hadamard manifold and
let $\om$ be a bounded domain with smooth boundary $\partial\Omega$
in $Q^n$. Let $\beta\geq 0$ and $\rho$ be a continuous positive
function on $\pa\om$. Then the first $n$ nonzero  eigenvalues of the
eigenvalue problem \be\label{th2.1}\left\{\begin{array}{l}
\ov{\Delta} u =0 \ \ \ \ \ \
\ \ \ \ \ \ \ \ \ \ \ \ \ \ \ \ \ \ \ {\rm in  \ }\ \om,\\
-\beta\Delta u+\pa_{\nu} u= \laa \rho u\ \ \ \  \qquad {\rm on  \ }
\pa \om
\end{array}\right.
\en satisfy \be\label{th2.2}
\sum_{i=1}^n\frac{1}{\lambda_{i,\b}}+\sum_{i=1}^{n-1}\frac{\b|\p\Omega|}{\lambda_{i,\b}|\om|}
\geq\fr{n^2|\om|}{\int_{\pa\om} \rho^{-1}}. \en
Moreover, when
$\rho=\rho_0$ is constant, equality in (\ref{th2.2}) holds  if and
only if $\om$ is isometric to an $n$-dimensional Euclidean ball.
\end{theorem}

Using the monotonicity of eigenvalues $\lambda_{i,\beta}$'s and
Theorem \ref{th2} immediately, we can obtain:

\begin{corollary}
Under the assumptions of Theorem \ref{th2}, we have
\begin{eqnarray} \label{add-5}
\lambda_{1,\beta}(\Omega)\leq\left[1+\beta\frac{(n-1)|\partial\Omega|}{n|\Omega|}\right]\cdot\frac{\int_{\pa\om}
\rho^{-1}}{n|\Omega|}.
\end{eqnarray}
Moreover, when $\rho=\rho_0$ is constant, equality in (\ref{add-5})
holds if and only if  $\om$ is isometric to an $n$-dimensional
Euclidean ball.
\end{corollary}

In the second part of the present paper, we consider eigenvalues of
the Laplacian acting on functions on closed hypersurfaces in a
Euclidean space or a Hadamard manifold, and can prove the following
two conclusions.

\begin{theorem} \label{th3}
Let $M$ be a connected closed hypersurface in $\R^{n}$ with $n\geq
3$.

(i) The first $(n-1)$ nonzero eigenvalues of the Laplacian on $M$
satisfy
\begin{eqnarray}\label{th3.1}
\sum_{j=1}^{n-1}\la_j^{\frac{1}{2}}\leq
(n-1)^{\frac{3}{2}}\left(\fr{\int_M
|{\H}|^2}{|M|}\right)^\frac{1}{2},
\end{eqnarray}
where ${\H}$ denotes the mean curvature vector of $M$ in $\R^{n}$
and, as before, $|M|$ is the area of $M$. The equality holds in
(\ref{th3.1}) if and only if $M$ is a hypersphere.

(ii) If $M$ is embedded and bounds a region $\om$, then
\be\label{the3.2} \sum_{i=1}^{n-1}
\fr{1}{\la_i}\geq\fr{n\omega_n}{|M|}\left(\fr{|\om|}{\omega_{n}}\right)^{1+\fr
1n}\left[1+\fr{(n+1)\left(2^{1/n}-1\right)}{4n}\left(\fr{|\om\Delta
B|}{|B|}\right)^2\right],  \en where $B$ is the ball of volume
$|\om|$ and with the same center of mass than $M$. Equality holds in
(\ref{the3.2}) if and only if $M$ is a hypersphere.
\end{theorem}

\begin{theorem} \label{th4}
  Let  $Q^n$ be a  Hadamard $n$-manifold ($n\geq 3$). Let $M$ be a connected closed hypersurface
embedded   in $Q^n$ and  $\om$ be the region bounded by $M$. Then
the first $(n-1)$ nonzero eigenvalues  of the Laplacian of $M$
satisfy\be\label{the4.1} \sum_{i=1}^{n-1}
\fr{1}{\la_i}\geq\fr{n^2|\om|^2}{|M|^2}. \en Equality holds in
(\ref{the4.1}) if and only if $\om$ is isometric to an
$n$-dimensional Euclidean ball.
\end{theorem}

\begin{remark}
\rm{ (1) By applying Theorem \ref{th3} and the monotonicity of
eigenvalues $\lambda_{i}$'s, we separately have:\footnote{Of course,
it should be ``\emph{under the assumptions of Theorem \ref{th3}}".}
\begin{eqnarray}
&&(i)\qquad \lambda_{1}(M)\leq(n-1)\fr{\int_M
|{\H}|^2}{|M|}; \label{add-6}\\
&&(ii)\qquad
\lambda_{1}(M)\leq\frac{(n-1)|M|}{n|\Omega|}\left(\frac{w_{n}}{|\Omega|}\right)^{\frac{1}{n}}
\left[1+\fr{(n+1)(2^{1/n}-1)}{4n}\left(\fr{|\om\Delta
B|}{|B|}\right)^2\right]^{-1}, \nonumber
\end{eqnarray}
and moreover, the equality holds implies the rigidity described as
in Theorem \ref{th3}. The eigenvalue estimate (\ref{add-6}) is
actually the well-known Reilly's inequality. In fact, estimate
(\ref{add-6}), together with the corresponding rigidity, can be
extended to the case of codimension $\geq2$, and this is actually
the main result of the influential paper \cite{RR}. Besides, there
is one more thing we prefer to mention here, that is, Mao and his
collaborators \cite[Theorem 1.11]{dmwxz} successfully gave a sharp
lower bound for the sum of the reciprocals of the first $n$ nonzero
eigenvalues of the Laplacian on a closed $n$-submanifold immersed in
a Euclidean space, and then Reilly's inequality follows as a direct
consequence. \\
 (2) By applying Theorem \ref{th4} and the monotonicity of
eigenvalues $\lambda_{i}$'s, we get\footnote{Of course, it should be
``\emph{under the assumptions of Theorem \ref{th4}}".}
\begin{eqnarray*}
\lambda_{1}(M)\leq\frac{(n-1)|M|^2}{n^{2}|\Omega|^2},
\end{eqnarray*}
 and equality holds implies the rigidity described as
in Theorem \ref{th4}.
 }
\end{remark}

Let $\tau>0$ be a positive constant, and let $\om$ be a bounded
domain in $\R^{n}$ with boundary $\partial\Omega$. As before, let
$\overline{\nabla}$, $\overline{\Delta}$ the gradient and the
Laplace operators in $\mathbb{R}^{n}$, respectively. Consider the
following
 Steklov-type eigenvalue problem
\begin{eqnarray}\label{th6.0}
\left\{\begin{array}{ccc} \ov{\D}^2 u-\tau \ov{\D} u
=0&&~\mbox{in} ~~ \om, \\[2mm]
\frac{\p^2u}{\p\nu^2}=0 &&~~\mbox{on}~~\partial \om,\\[2mm]
\tau\frac{\p u}{\p \nu}-\mathrm{div}_{\p \Omega}\(\ov{\nabla}^2
u\cdot\nu\)_{\partial\Omega}-\frac{\p\ov{\D} u}{\p\nu}=\lambda u
&&~~\mbox{on}~~\partial \om,
\end{array}\right.
\end{eqnarray}
where  $\div_{\p \Omega}$ denotes the tangential divergence operator
on $\p \Omega$, $\ov{\nabla}^2 u$ is the Hessian of $u$,
$(\overline{\nabla}^{2}u\cdot\nu)_{\partial\Omega}$ stands for the
projection of $\overline{\nabla}^{2}u\cdot\nu$ to the tangent bundle
of $\partial\Omega$. This problem (\ref{th6.0}) has a discrete
spectrum whose elements (i.e., eigenvalues) can be listed
non-decreasingly as follows
$$0=\lambda_{0,\tau}<\lambda_{1,\tau}\leq\lambda_{2,\tau}\leq\cdots\leq\lambda_{k,\tau}\leq\cdots\rightarrow\infty.$$
The eigenvalue $0$ is simple and its eigenfunctions are nonzero
constant functions. By the variational principal, the $k$-th $(k\geq
1)$ eigenvalue $\lambda_{k,\tau}$ can be characterized as follows
 \begin{eqnarray}\label{th6.1}
\lambda_{k,\tau}=\mathrm{min}\Bigg{\{}\frac{\int_{\om}\(|\ov{\na}^2u|^2+\tau|\overline{\na}
u|^2\)}{\int_{\pa\om} u^2}\Bigg| u\in H^2(\om), u\neq0, \nonumber\\
\int_{\pa\om} u u_j=0, j=0,\cdots,k-1\Bigg{\}},
\end{eqnarray}
where $u_{j}$ is the eigenvalue function of $\lambda_{j,\tau}$. We
can prove the following sharp upper bound estimate:

\begin{theorem} \label{th6} For
the eigenvalue problem (\ref{th6.0}), we have
\begin{eqnarray*}\label{th6.2}
\frac{1}{n-1}\sum_{j=1}^{n-1}\lambda_{j,\tau}^{\frac{1}{2}}\leq
\frac{\left\{ n\tau |\Omega|\inpo |\H|^2\right\}^\frac{1}{2}}{|\p
\Omega|},
\end{eqnarray*}
where $\H$ is the mean curvature vector of $\p \Omega$ in $\R^n$.
Equality in (\ref{th6.2}) holds  if and only if $\Omega$ is a ball.
\end{theorem}

\begin{remark}
\rm{ (1) In \cite{LM}, the authors therein used the operator
$\mathrm{Proj}_{\partial\Omega}\left[(\overline{\nabla}^{2}u)\nu\right]$
to denote the projection  of $(\overline{\nabla}^{2}u)\nu$ onto the
space tangent to $\partial\Omega$, which obviously has the same
meaning as
$(\overline{\nabla}^{2}u\cdot\nu)_{\partial\Omega}$ here. \\
 (2) By applying Theorem \ref{th6} and the monotonicity of
eigenvalues $\lambda_{i,\tau}$'s, one has
\begin{eqnarray*}
\lambda_{1,\tau}(\Omega)\leq \frac{n\tau |\Omega|\inpo |\H|^2}{|\p
\Omega|^2},
\end{eqnarray*}
with equality  if and only if $\Omega$ is a ball. \\
 (3) The Steklov-type eigenvalue problem (\ref{th6.0}) has been
 studied by Buoso and Provenzano \cite{bp} already, and when $n=2$, it can be used
 to describe transverse vibrations of a thin plate, which arises in the theory of linear elasticity (see \cite[Section 2]{bp} for
 details). One can also check \cite[Section 3]{bp} for a detailed
 explanation about the existence of discrete spectrum and the
 characterization (\ref{th6.1}) of the $k$-th eigenvalue of the eigenvalue problem
 (\ref{th6.0}). Besides, \cite[Corollary 5.20]{bp} tells us that:
\begin{itemize}
\item \emph{Among all bounded domains of class $C^1$ with fixed measure, the ball maximizes
the first nonnegative eigenvalue of problem (\ref{th6.0}), that is,
$\lambda_{1,\tau}(\Omega)\leq\lambda_{1,\tau}(B)$ with $B$ is the
Euclidean ball having the same measure as $\Omega$.}
\end{itemize}
Very recently, this interesting spectral isoperimetric inequality
has been improved by Li and Mao -- see \cite[Theorem 1.1]{LM1}.
 }
\end{remark}

\section{Proofs of Theorems \ref{th1} and \ref{th2}}
\renewcommand{\thesection}{\arabic{section}}
\renewcommand{\theequation}{\thesection.\arabic{equation}}
\setcounter{equation}{0} \label{intro}

In this section, we will give the proofs of Theorems \ref{th1} and
\ref{th2}. First, we recall the following lemma from \cite{bpr}.

\begin{lemma}\label{le1} Let $\om$ be a bounded, open Lipschitz set
in $\R^n$. Then \be\label{sec2.1} \int_{\pa\om} |x|^2  \geq R^2|\pa
B_R|\(1+\fr{(n+1)(2^{1/n}-1)}{4n}\left(\fr{|\om\Delta
B_R|}{|B_R|}\right)^2\), \en where $B_R$ is the ball centered at the
origin such that $|B_R|=|\om|$.
\end{lemma}

$\\${\it Proof of Theorem \ref{th1}.} Let $u_0, u_1, u_2,\cdots$ be
the orthonormal (in $L^2(\pa \om)$ sense) eigenfunctions
corresponding to the eigenvalues
$0=\lambda_{0,\b}<\lambda_{1,\b}\leq \lambda_{2,\b}\leq \cdots $ of
the eigenvalue problem (\ref{th1.1}), that is,
\begin{eqnarray*}\label{pth1.1}\left\{\begin{array}{l} \ov{\Delta} u_i =0 \ \ \
\
\ \ \ \ \ \ \ \ \ \ \ \ \ \ \ \ \ \ \ \ {\rm in  \ }\ \om,\\
-\beta\Delta
u+\pa_{\nu} u_i= \laa_i u_i\ \ \ \ \quad  {\rm on  \ } \pa \om,\\
\int_{\pa \om} u_iu_j=\delta_{ij}.
\end{array}\right.
\end{eqnarray*}
 By (\ref{vc}), the eigenvalues $\lambda_{i,\b}, i=1, 2,\cdots,$ are
characterized by \be \label{pth1.2}
\lambda_{i,\b}=\underset{\underset{u\perp{\rm span}\{u_0,\cdots,
u_{i-1}\}}{u\in H(\om)\setminus\{0\}}}{\min}
\fr{\int_{\om}|\ov{\na}u|^2 +\beta\int_{\pa\om} |\na
u|^2}{\int_{\pa\om} u^2}.\en We need to choose nice trial functions
$\phi_i$ for each of the eigenfunctions $u_i$ and insure that these
are orthogonal to the preceding eigenfunctions $u_0,\cdots,
u_{i-1}$. For the $n$ trial functions $\phi_1, \phi_2, \cdots,
\phi_n,$ we simply choose the $n$ coordinate functions
 \be
\phi_i=x_i, \ \ {\rm for}\ \ i=1,\cdots, n, \nonumber\en
 but before
we can use these we need to make adjustments  so that $\phi_i
\perp{\rm span}\{u_0,\cdots, u_{i-1}\}$ in $L^2(\pa\om)$. Simply, by
translating the origin appropriately we can assume that \be
\int_{\pa\om} x_i=0, \ i=1,\cdots,n, \nonumber \en
 that is, $x_i\perp u_0$
(which is actually just the constant function $1/\sqrt{|\pa\om|}$).
Nextly we show that a suitable rotation of axes can be made so as to
insure that \be\label{pth1.} \int_{\pa\om}\phi_j u_i=\int_{\pa\om}
x_j u_i=0, \en for $j=2,3,\cdots, n$ and $i=1,\cdots,j-1$. To see
this, we define an $n \times n$ matrix $Q=\(q_{ji}\)_{n\times n},$
where $q_{ji}=\int_{\p\om} x_j u_i$, for $i,j=1,2,\cdots,n.$ Using
the orthogonalization of Gram and Schmidt (QR-factorization
theorem), we know that there exist an upper triangle matrix
$T=(T_{ji})_{n\times n}$ and an orthogonal matrix
$U=(a_{ji})_{n\times n}$ such that $T=UQ$, i.e.,
\begin{eqnarray*}
T_{ji}=\sum_{k=1}^n a_{jk}q_{ki}=\inpo \sum_{k=1}^n  a_{jk}x_k  u_i
=0,\ \  1\leq i<j\leq n.
\end{eqnarray*}
Letting $y_j=\sum_{k=1}^n  a_{jk}x_k$, we get
\begin{eqnarray*}\label{c2}
\inpo y_j  u_i =\inpo \sum_{k=1}^n a_{jk}x_k u_i
=0,\ \  1\leq i<j\leq n.
\end{eqnarray*}
Since $U$ is an orthogonal matrix,  $y_1, y_2, \cdots, y_n$ are also
coordinate functions on $\R^n$. Therefore, denoting these coordinate
functions still by $x_1, x_2,\cdots, x_n$, we arrive at the
condition (\ref{pth1.}).

It then follows from (\ref{pth1.2}) that for each fixed
$i=1,\cdots,n$,
\begin{eqnarray}\label{d5}
\lambda_{i,\b}\inpo x_i^2 \leq\ino |\ov{\n} x_i|^2 +\beta\inpo|\n
x_i|^2= |\Omega|+\beta\inpo|\n x_i|^2, \nonumber
\end{eqnarray}
with equality holding if and only if \be\label{pth}\beta\Delta x_i
+\pa_{\nu} x_i = -\laa_{i, \beta} x_i  \ \ \ {\rm on}\ \ \ \pa\Om.
\en
 Hence,
\begin{eqnarray}\label{pth1.3}
\inpo
x_i^2\leq\frac{|\Omega|}{\lambda_{i,\b}}+\frac{\beta}{\lambda_{i,\b}}\inpo|\n
x_i|^{2}.
\end{eqnarray}
Observing
\begin{eqnarray*}
\sum_{i=1}^n|\n x_i|^2=n-1,\ |\n x_i|^2\leq1,
\end{eqnarray*}
we get
\begin{eqnarray}\label{pth1.4}
\nonumber\sum_{i=1}^n\frac{|\n
x_i|^2}{\lambda_{i,\b}}&=&\sum_{i=1}^{n-1}\frac{|\n
x_i|^2}{\lambda_{i,\b}}+\frac{|\n x_n|^2}{\lambda_{n,\b}}
\\\nonumber&=&\sum_{i=1}^{n-1}\frac{|\n x_i|^2}{\lambda_{i,\b}}+\frac{1}{\lambda_{n,\b}}\sum_{i=1}^{n-1}\(1-|\n x_i|^2\)
\\\nonumber&\leq&\sum_{i=1}^{n-1}\frac{|\n x_i|^2}{\lambda_{i,\b}}+\sum_{i=1}^{n-1}\frac{1}{\lambda_{i,\b}}\(1-|\n x_i|^2\)
\\&=&\sum_{i=1}^{n-1}\frac{1}{\lambda_{i,\b}}.
\end{eqnarray}
Combining (\ref{pth1.3}) and (\ref{pth1.4}), we have
\begin{eqnarray}\label{c8}
\sum_{i=1}^n\inpo
x_i^2\leq\sum_{i=1}^n\frac{|\Omega|}{\lambda_{i,\b}}+\sum_{i=1}^{n-1}\frac{\b|\p\Omega|}{\lambda_{i,\b}}.
\end{eqnarray}
Substituting (\ref{le1}) into (\ref{c8}), one gets (\ref{th1.2}). If
the equality holds in (\ref{th1.2}), then the inequality
(\ref{pth1.4}) must take equality sign and (\ref{pth}) holds. It is
easy to see from the equality case of (\ref{pth1.4}) at any point of
$\pa\om$ that \be\lambda_{1,\b}=\lambda_{2,\b}\cdots=\lambda_{n,\b}.
\nonumber \en
 It then follows that the
position vector $x=(x_1,\cdots, x_n)$ when restricted on $\pa\Om$
satisfies \be\label{p1} \Delta x=(\Delta x_1,\cdots, \Delta x_n) =
-\fr{1}{\beta}\nu-\fr{\laa_{1, \beta}}{\beta}(x_1,\cdots, x_n). \en
On the other hand,  it is well known that \be\label{p2} \Delta x
=(n-1){\mathbf H}, \en where ${\mathbf H}$ is the mean curvature
vector of $\pa \Om$ in ${\R}^n$. Combining (\ref{p1}) and
(\ref{p2}), we have \be\label{p3} x=-\fr
1{\laa_{1,\beta}}\nu-\fr{(n-1)\beta}{\laa_{1,\beta}}{\mathbf H} \ \
\ \ {\rm on\ \ \ }\pa \Om. \en Consider the function $g = |x|^2 :
M\rightarrow {\R}$. It is easy to see from (\ref{p3}) that \be\no Z
g =2\langle Z, x\rangle = 0, \ \ \ \forall Z\in
{\mathfrak{X}}(\pa\Om), \en
 where ${\mathfrak{X}}(\pa\Om)$ is the tangent bundle of
 $\partial\Omega$.
 Thus $g$ is a constant function and so
$\pa \Om$ is a hypersphere. Theorem \ref{th1} follows.  \hfill
$\Box$

$\\${\it Proof of Theorem \ref{th2}.} Let $\{u_i\}_{i=0}^{+\infty}$
be an orthonormal set of eigenfunctions corresponding to the
eigenvalues $\{\laa_i\}_{i=0}^{+\infty}$
 of the eigenvalue problem (\ref{th2.1}),
that is, \begin{eqnarray*}\label{pth2.1}\left\{\begin{array}{l}
\ov{\Delta} u_i =0
\ \ \ \ \ \ \ \ \ \ \ \ \ \ \ \ \ \ \ \quad  \ {\rm in  \ }\ \om,\\
-\beta\Delta
u+\pa_{\nu} u_i= \laa_{i, \beta} u_i\ \ \ \ {\rm on  \ } \pa \om,\\
\int_{\pa \om} \rho u_iu_j=\delta_{ij}.
\end{array}\right.
\end{eqnarray*}
 By (\ref{vc}), we have  \be \label{pth2.2}
\lambda_{i,\b}=\underset{\underset{\int_{\pa \om} \rho uu_l=0, l=0,
1,...,i-1}{u\in H(\om)\setminus\{0\}}}{\min}
\fr{\int_{\om}|\ov{\na}u|^2 +\beta\int_{\pa\om} |\na
u|^2}{\int_{\pa\om} \rho u^2}.\en

The idea is to use globally defined coordinate functions suitably
chosen on $Q^n$ as trial functions for the first $n$ nonzero
eigenvalues of the problem (\ref{th2.1}). To do this, let us denote
by $\langle, \rangle$ the Riemannian metric on $Q^n$. For any $p\in
Q^n$, let ${\rm exp}_p$ and $UQ_{p}^n$ be the exponential map and
unit tangent space of $Q^n$ at $p$, respectively. Let $\{e_1,\cdots,
e_n\}$ be an orthonormal basis of $T_{p}Q^n$, the tangent space of
$Q^n$ at $p$, and $y : M\ri \R^n$ be the Riemannian normal
coordinates on $Q^n$ determined by $(p; e_1,\cdots,e_n)$. It follows
from  the Cartan-Hadamard theorem (see, e.g., \cite{car}) that $y$
is well-defined on all of $Q^n$ and is a diffeomorphism of $Q^n$
onto $\R^n$. We can choose $p$ and $\{e_1,\cdots, e_n\}$ so that the
respective coordinate functions $ y^i : M\ri \R, i=1,\cdots,n,$ of
$y : Q^n \ri \R^n$ satisfy \be \label{pth2.3}\int_{\pa \om} \rho y_i
=0. \en In fact, parallelly translate the frame $\{e_1,\cdots,
e_n\}$ along every geodesic emanating from $p$ and thereby obtain a
differentiable orthonormal frame field $\{E_1, \cdots, E_n\}$ on
$Q^n$. Let $y^q : M\ri \mathbb{R}^n$ denote the Riemannian normal
coordinates of $Q^n$ determined by $\{E_1, \cdots, E_n\}$ at $q$,
and let $y^q_i, i=1,\cdots, n,$ be the coordinate functions of
$y^q$. By definition, $ y^q_i : M\ri \mathbb{R}$ is given by
$y_i^q(z)=\langle {\rm exp}_q^{-1}(z), E_i(q)\ra$, $i=1,\cdots,n.$
Then
 \be Y(q)=\sum_{i=1}^n \left\{\int_{\pa\om}\rho y^q_i \right\}
E_i(q) \nonumber \en
 is a continuous vector field on $Q^n$. If we
restrict $Y$ to a geodesic ball $B$ containing $\overline{\om}$ then
the convexity of $B$ implies that on the boundary of $B$, $Y$ points
into $B$. The Brouwer fixed point theorem  then implies that $Y$ has
a zero. So we may assume that $p$ and $\{e_1,\cdots, e_n\}$ actually
satisfy (\ref{pth2.3}).

For any $w\in UQ^n_p$, let $\theta_w$ be the function on $Q^n$
defined by $\theta_w(z)=\langle {\rm exp}_p^{-1}(z), w\ra $. Then
(\ref{pth2.3}) is equivalent to say that for $i=1, \cdots, n,$ \be
\int_{\pa\om} \rho\theta_{e_i} =0. \nonumber \en Thus for any
$\sigma\in UQ^{n}_p$, \be\label{pth2..} \int_{\pa\om}
\rho\theta_{\sigma}=0. \en Next, we show that there exists an
orthonormal basis $\{\sigma_j\}_{j=1}^n$ of $T_pQ^n$ such that the
coordinate functions $x_{\a}=\theta_{\sigma_{\alpha}},
\a=1,\cdots,n$ of the Riemannian normal coordinate system of $Q^{n}$
determined by $(p; \sigma_1,\cdots, \sigma_n)$ satisfy
\be\label{pth2.5} \int_{\pa\om} \rho x_i u_j=0, \ \ {\rm for} \ \
i=2, 3,\cdots, n,\ {\rm and}\ j=1,\cdots, i-1.\en In fact, let us
consider the $n\times n$ matrix $P=\(p_{\a\b}\)_{n\times n}$, where
$p_{\a\b}=\int_{\pa\om} \rho y_\a u_\b=\int_{\pa\om}
\rho\theta_{e_{\alpha}} u_\b$, for $\a,\b=1,2,\cdots,n.$ Using the
same discussion as that in the proof of Theorem \ref{th1}, we can
find
 an orthogonal matrix
$U=(a_{\a\b})_{n\times n}$ such that
\begin{eqnarray*}
\inpo \sum_{\gamma=1}^n a_{\a\gamma}\rho y_\gamma  u_\b =0,\ \
1\leq\b<\a\leq n.
\end{eqnarray*}
Setting $\sigma_\a=\sum_{\gamma=1}^n  a_{\a\gamma}e_\gamma$, we know
that  $\{\sigma_j\}_{j=1}^n$ is an orthonormal basis of $T_{p}Q^{n}$
and that the condition (\ref{pth2.5}) is satisfied. Therefore, we
have from (\ref{pth2.2}), (\ref{pth2..}) and (\ref{pth2.5}) that
\be\label{pth2.7} \lambda_{i,\b}\int_{\pa\om} \rho x_i^2\leq
\int_{\om}|\ov{\na}x_i|^2 +\beta\int_{\pa\om} |\na x_i|^2, \quad
i=1,\cdots,n. \en Let $\left\{\fr{\pa}{\pa x_k}, k=1,\cdots,
n\right\}$ be the natural basis of the tangent spaces associated
with the coordinate chart $x$ and let $g_{kl}= \langle\fr{\pa}{\pa
x_k}, \fr{\pa}{\pa x_l}\rangle $, $k, l=1,\cdots, n$. Since $Q^n$
has non-positive sectional curvature, the Rauch comparison theorem
(see, e.g., \cite{car}) implies that the eigenvalues of the matrix
$(g_{kl})_{n\times n}$ are all $\geq 1$. Thus the eigenvalues of
$(g^{kl})_{n\times n}=: (g_{kl})_{n\times n}^{-1}$ are $\leq 1$ and
so we have $g^{kk}\leq 1, \ k=1,\cdots, n$. Thus we have for
$i=1,\cdots, n,$ \be\label{pth2.8}|\na x_i|^2\leq |\ov{\na} x_i|^2 =
g^{ii} \leq 1. \en We claim that \be\label{pth2.9} \sum_{i=1}^n |\na
x_i|^2 \leq n-1 \ \ {\rm on}\ \ \pa\om. \en In fact, for a fixed
point $z\in \pa\om$, let $\{v_1,\cdots,v_{n-1}\}$ be an orthonormal
basis of $T_z(\pa\om)$. From the fact that the eigenvalues of
$(g_{ij})_{n\times n}$ are all $\geq 1$, we conclude that \be
\sum_{i=1}^n |\na x_i|^2&=&\sum_{i=1}^n \sum_{\alpha=1}^{n-1}
(v_{\alpha}x_i)^2 \no\\
\no &=& \sum_{\alpha=1}^{n-1}\sum_{i,j=1}^n (v_{\alpha}x_i)
\delta_{ij}(v_{\alpha}x_j)\\ \no &\leq&
\sum_{\alpha=1}^{n-1}\sum_{i,j=1}^n (v_{\alpha}x_i)
g_{ij}(v_{\alpha}x_j)\\ \no &=& \sum_{\alpha=1}^{n-1}
|v_{\alpha}|^2=n-1. \en Using (\ref{pth2.7})-(\ref{pth2.9}) and the
same arguments as in the proof of Theorem \ref{th1}, we obtain
\be\label{pth2.10} \int_{\pa\om} \rho r^2 = \int_{\pa\om} \rho
\sum_{i=1}^n x_i^2\leq
\sum_{i=1}^n\frac{|\om|}{\lambda_{i,\b}}+\sum_{i=1}^{n-1}\frac{\b|\p\Omega|}{\lambda_{i,\b}},
\en where $r = d(p, \cdot): \om\ri {\R}$ denotes the distance
function from $p$. The Cauchy-Schwarz inequality implies that
\be\label{pth2.11} \int_{\pa\om} \rho r^2 \geq\fr{(\int_{\pa \om} r
)^2}{\int_{\pa\om} \rho^{-1}} \en with equality holding if and only
if $r\rho=const.$ on $\pa\om$. From the Laplace comparison theorem
(cf. \cite{sy}), we have  \be\label{pth2.12} \D r^2 \geq 2n. \en
Integrating (\ref{pth2.12}) on $\om$ and using the divergence
theorem, we get \be\label{pth2.13} n|\om|\leq \int_{\pa \om}
r\langle \na r, \nu\ra\leq \int_{\pa \om} r |\na r|=\int_{\pa \om}
r. \en Combining (\ref{pth2.10}), (\ref{pth2.11}) and
(\ref{pth2.13}), we get (\ref{th2.2}).

Assume now that $\rho=\rho_0$ is constant and that  the equality
holds in (\ref{th2.2}). In this case, we  have $r|_{\pa\om}=const.$
and so $\om$ is a geodesic ball with center $p$. Also, we have \be
\Delta r^2|_{\om}= 2n. \en It then follows from the equality case in
the Laplace comparison theorem (cf. \cite{sy}) and the Cartan's
theorem (see, e.g., \cite{car}) that $\om$ is isometric to an
$n$-dimensional Euclidean ball. This completes the proof of Theorem
\ref{th2}.\hfill $\Box$

\section{Proofs of Theorems \ref{th3}-\ref{th6}}
\renewcommand{\thesection}{\arabic{section}}
\renewcommand{\theequation}{\thesection.\arabic{equation}}
\setcounter{equation}{0} \label{intro}

In this section, we will give the proofs of Theorems
\ref{th3}-\ref{th6} in detail.

$\\${\it Proof of Theorem \ref{th3}.} Let $\{u_i\}_{i=0}^{+\infty}$
be the orthonormal system of eigenfunctions  corresponding to the
eigenvalues \be\label{pth3.1} 0=\la_0<\la_1\leq
\la_2\leq\cdots\ri\infty \nonumber \en of the Laplacian of $M$, that
is, \be\label{th3.2} \Delta u_i = -\la_i u_i, \quad  \int_M u_iu_j
=\delta_{ij}. \en We have $u_0=1/\sqrt{|M|}$ and for each
$i=1,\cdots,$ \be \label{pth3.3} \la_i=\underset{u\neq 0, \int_M
uu_j=0, j=0,\cdots i-1}{\min}\fr{\int_M |\na u|^2}{\int_M u^2}. \en
Let $x_1, \cdots, x_{n}$ be the coordinate functions on $\R^{n}$. By
using the same arguments as in the proof of Theorem \ref{th1}, we
can assume that
\begin{eqnarray*}
\ino x_i u_j =0,\  i=1,\cdots,n;\ j=0,\cdots,i-1
\end{eqnarray*}
and so we have \be\label{pth3.4} \la_i\int_M  x_i^2\leq \int_M |\na
x_i|^2, \ i=1,\cdots, n,
\end{eqnarray}
with equality holding if and only if \be\label{pth3.5} \D x_i=
-\la_i x_i. \nonumber \en Since \be\label{pth3.51} |\na x_i|^2\leq
1, \quad \sum_{i=1}^n |\na x_i|^2=n-1, \en we have\be \label{pth3.6}
\no \sum_{i=1}^n \la_i^{\fr 12}|\na x_i|^2&=&\sum_{i=1}^{n-1}
\la_i^{\fr 12}|\na x_i|^2+\la_n^{\fr 12}|\na x_n|^2\\ \no &=&
\sum_{i=1}^{n-1} \la_i^{\fr 12}|\na x_i|^2+\la_n^{\fr
12}\sum_{i=1}^{n-1}(1-|\na x_i|^2)
\\ \no &\geq & \sum_{i=1}^{n-1}
\la_i^{\fr 12}|\na x_i|^2+\sum_{i=1}^{n-1}\la_i^{\fr 12}(1-|\na
x_i|^2)
\\ \no &=& \sum_{i=1}^{n-1}\la_i^{\fr 12},
\en which gives \be \label{pth3.6}\sum_{i=1}^n \la_i^{\fr 12}\int_M
|\na x_i|^2 \geq |M|\sum_{i=1}^{n-1}\la_i^{\fr 12}. \en For any
positive constant $\delta$, we have from Schwarz inequality and
(\ref{pth3.4}) that \be\label{pth3.7} \no \sum_{i=1}^n \la_i^{\fr
12}\int_M |\na x_i|^2 &=&\sum_{i=1}^n \la_i^{\fr 12}\int_M (-x_i \D
x_i)\\ \no &\leq& \fr 12\sum_{i=1}^n\left\{\delta\la_i\int_M x_i^2 +
\fr 1{\delta}\int_M (\D x_i)^2\right\}\\  &\leq& \fr 12\left\{
(n-1)\delta |M|+ \fr 1{\delta}\int_M (n-1)^2 |{\bf H}|^2\right\}.
 \en
Taking $$\delta =\left\{\fr{\int_M (n-1)|{\bf
H}|^2}{|M|}\right\}^{\fr 12}$$ in (\ref{pth3.7}) and  using
(\ref{pth3.6}), we get (\ref{th3.1}) with equality holding if and
only if $M$ is a hypersphere.

On the other hand, one can use (\ref{pth3.4}), (\ref{pth3.51}) and a
similar argument as that in the proof of Theorem \ref{th1} to obtain
 \be
\sum_{i=1}^n \int_M x_i^2 \leq |M|\sum_{i=1}^{n-1} \fr 1{\la_i},
 \en
which, combining with (\ref{sec2.1}) gives (\ref{the3.2}). Also, one
can deduce as in the proof of the final part of Theorem \ref{th1}
that  equality holds in (\ref{th3.1}) if and only if $M$ is a sphere
of $\R^{n}$. This completes the proof of Theorem \ref{th3}. \hfill
$\Box$

$\\${\it Proof of Theorem \ref{th4}.} From the proof of Theorem
\ref{th2}, we know that there exist a point $p\in Q^n$ and an
orthonormal basis $\{e_i\}_{i=1}^n$ of $T_pQ^n$ such that the
coordinate functions of the Riemannian normal coordinate system
determined by $\{p; (e_1,\cdots, e_n)\}$ satisfy \be \label{pth4.1}
\int_M x_i \phi_j =0, \ i=1,\cdots,n, \ j=0,\cdots, i-1, \no \en
where $\{\phi_j\}_{j=0}^{+\infty}$ are orthonormal eigenfunctions
corresponding to the eigenvalues \be 0=\la_0<\la_1\leq \la_2\leq
\cdots \ri +\infty \no\en of the Laplacian of $M$. Thus we have \be
\la_i\int_M x_i^2 \leq \int_M |\na x_i|^2, \ i=1,...,n. \no\en Using
a similar argument as that in the proof of Theorem \ref{th2}, we
conclude \be\sum_{i=1}^{n-1} \fr 1{\la_i}\geq \fr{\int_M r^2}{|M|}
\geq \fr{n^2|\om|^2}{|M|^2}, \no \en with equality holding if and
only if $\om$ is isometric a ball in $\R^n$, where $r=d(p,\cdot):
\om\ri \R$ is the distance function from the point $p$. The proof of
Theorem \ref{th4} is finished.\hfill $\Box$

$\\$ {\it Proof of Theorem \ref{th6}.}  Let $x=(x_1, \cdots, x_n)$
be the coordinate functions on $\R^n$. Since $\Omega$ is a bounded
domain with smooth boundary $\p \Omega$ in $\R^n$, we can regard $\p
\Omega$ as a closed hypersurface
 of $\R^n$ without boundary.
Let $\D,\n$ be the Laplace operator and the gradient operator on $\p
\om$, respectively. The position vector $x=(x_1, \cdots, x_n)$ when
restricted on $\pa\Omega$ satisfies
\begin{eqnarray*}\label{c1}
\sum_{\a=1}^n \(\D x_\a\)^2=(n-1)^2|\H|^2,
\end{eqnarray*}
where $|\H|$ is the mean curvature of $\p \Omega$ in $\R^n$.

Let $u_i$ be the eigenfunction corresponding to eigenvalue
$\lambda_{i,\tau}$  such that $\{u_i\}_ {i=0}^{\infty}$ becomes an
orthonormal basis of $L^2(\partial \Omega)$, that is,
\begin{eqnarray*}
\left\{\begin{array}{ccc} \L^2 u_i-\tau\L u_i
=0,&&~\mbox{in}~~\Omega \\[2mm]
\frac{\p^2u_i}{\p\nu^2}=0 &&~~\mbox{on}~~\partial \Omega,\\[2mm]
\tau\frac{\p u_i}{\p \nu}-\mathrm{div}_{\p \Omega}\(\overline{\nabla}^2 u_i\cdot\nu\)_{\partial\Omega}-\frac{\p\L u_i}{\p\nu}=\lambda_{i,\tau} u_i&&~~\mbox{on}~~\partial \Omega,\\[2mm]
\inpo u_i u_j =\delta_{ij}.
\end{array}\right.
\end{eqnarray*}
Observe that $u_0= 1/\sqrt{|\pa \om|}$. We can assume as before that
\be \int_{\pa\om} x_i u_j=0, \ i=1,\cdots, n, \ j=0,\cdots, i-1.
\no\en It follows from the variational characterization
(\ref{th6.1}) that
\begin{eqnarray*}\label{c4}
\lambda_{i,\tau}\int_{\pa\om} x_i^2\leq \int_{\om}
\(|\ov{\na}^2x_i|^2 +\tau|\ov{\na} x_i|^2\) = \tau |\om|.
\end{eqnarray*}
Thus we have
\begin{eqnarray*}\label{c7}
\sum_{i=1}^n\lambda_{i,\tau}\int_{\pa\om} x_i^2\leq n\tau|\om|,
\end{eqnarray*}
which implies that
\begin{eqnarray}\label{pth5.4}
\no\sum_{i=1}^n\lambda_{i,\tau}^{\fr 12}\int_{\pa\om}|\na x_i|^2
&=&\sum_{i=1}^n\lambda_{i,\tau}\int_{\pa\om}(-x_i \D x_i)\\ \no
&\leq& \fr 12\sum_{i=1}^n\left\{\delta \lambda_{i,\tau}\int_{\pa\om}
x_i^2 +\fr 1{\delta} \int_{\pa\om}(\D x_i)^2\right\}\\  &\leq& \fr
12 \left\{ \delta n\tau|\om|+ \fr 1{\delta}
\int_{\pa\om}(n-1)^2|{\bf H}|^2\right\} \en for any $\delta >0.$

Using a similar argument as that in the proof of Theorem \ref{th3},
we know that \be\label{pth5.5} \sum_{i=1}^n \la_{i,\tau}^{\fr
12}\int_{\pa\om} |\na x_i|^2\geq |\pa\om|\sum_{i=1}^{n-1}
\la_{i,\tau}^{\fr 12}.\en Combining (\ref{pth5.4}) with
(\ref{pth5.5}) and taking
$$\delta=\left\{\frac{\inpo (n-1)^2 |\H|^2}{n\tau
|\Omega|}\right\}^{\frac{1}{2}},$$ we obtain
\begin{eqnarray*}\label{pth6.6}
\frac{1}{n-1}\sum_{j=1}^{n-1}\lambda_{j,\tau}^{\frac{1}{2}}\leq
\frac{\left\{ n\tau |\Omega|\inpo |\H|^2\right\}^\frac{1}{2}}{|\p
\Omega|}.
\end{eqnarray*}
If the equality holds in above inequality, we know that
\begin{eqnarray*}
\lambda_{1,\tau}=\cdots=\lambda_{n,\tau}\equiv \Lambda
\end{eqnarray*}
 and
\begin{eqnarray*}\label{c12}
\D x_i =-\delta\sqrt{\lambda_{n,\tau}} x_i, \ i=1,\cdots, n, \ \
{\rm on}\ \ \pa\om.
\end{eqnarray*}
  Since $\pa\Omega$ is a closed hypersurface  of $\R^n$, we conclude that   $\p \om$ is a round sphere.
 This completes the proof of Theorem \ref{th6}. \hfill $\Box$

\section*{Acknowledgements}
 F. Du was
supported by Hubei Key Laboratory of Applied Mathematics (Hubei
University), Research Team Project of Jingchu University of
Technology (Grant No. TD202006) and Research Project of Jingchu
University of Technology (Grant Nos. YB202010, ZX202002, ZX202006).
J. Mao was partially supported by the NSF of China (Grant No.
11801496), the Fok Ying-Tung Education Foundation (China), and Hubei
Key Laboratory of Applied Mathematics (Hubei University). Q. Wang
was partially supported by CNPq, Brazil (Grant No. 307089/2014-2).
C. Xia was partially supported by CNPq, Brazil (Grant No.
306146/2014-2).


\begin{thebibliography}{9999}
\bibitem{B} F. Brock, \emph{An isoperimetric inequality for eigenvalues of the Steklov problem}, ZAMM {\bf 81} (2001) 69--71.

\bibitem{bb} M. F. Betta, F. Brock, A. Mercaldo, M. R. Posteraro,  \emph{A weighted isoperimetric inequality and applications to
symmetrization}, J. Inequal. Appl. {\bf 4} (1999) 215--240.

\bibitem{bpr} L. Brasco, G. De Philippis, B.  Ruffini,  \emph{Spectral optimization for
the Stekloff-Laplacian: the stability issue}, J. Funct. Anal. {\bf
262} (2012) 4675--4710.

\bibitem{bp} D. Buoso, L. Provenzano,  \emph{A few shape optimization results for a biharmonic Steklov problem}, J. Differential Equations {\bf 259} (2015)  1778-1818.

\bibitem{ceg} B.  Colbois, A. El Soufi, A. Girouard, \emph{Isoperimetric control of
the Steklov spectrum}, J. Funct. Anal. {\bf 261} (2011) 1384--1399.


\bibitem{DKL} M. Dambrine, D. Kateb, J. Lamboley, \emph{An extremal eigenvalue problem for the Wentzell-Laplace operator},
 Ann. I. H. Poincar\'e Non Linear Analysis {\bf 33}(2) (2014) 409--450.

\bibitem{car} M. P. do Carmo, \emph{Riemannian Geometry}, Birkh\"{a}user, Boston, 1993.

\bibitem{DWX}F. Du, Q. Wang, C. Xia, \emph{Estimates for eigenvalues of the Wentzell-Laplace operator},
 J. Geom. Phys. {\bf 129} (2018) 25--33.

\bibitem{dmwxz} F. Du, J. Mao, Q. Wang, C. Xia, Y. Zhao, \emph{Estimates for eigenvalues of the Neumann and Steklov
problems}, available online at arXiv:1902.08998v2.

\bibitem{E1} J. F. Escobar, \emph{The geometry of the first non-zero Steklov eigenvalue}, J. Funct. Anal. {\bf 150} (1997) 544--556.
.
\bibitem{E2} J. F. Escobar, \emph{An isoperimetric inequality and the first Steklov eigenvalue}, J. Funct. Anal. {\bf 165} (1999) 101--116.

\bibitem{FS} A. Fraser, R. Schoen, \emph{The first Steklov eigenvalue, conformal geometry, and minimal surfaces}, Adv. Math. {\bf 226} (2011) 4011--4030.

\bibitem{G} G. R. Goldstein, \emph{Derivation and physical interpretation of general boundary conditions}, Adv.
Differential Equations {\bf 11} (2006) 457--480.


\bibitem{HPS}J. Hersch, L. E. Payne,  M. M. Schiffer, \emph{Some inequalities for Stekloff eigenvalues}, Arch. Ration. Mech. Anal. {\bf 57} (1974) 99--114.

\bibitem{hx} G. N. Hile, Z. Xu, \emph{Inequalities for sums of the
reciprocals of eigenvalues}, J. Math. Anal. Appl. {\bf 180} (1993)
412--430.

\bibitem{kk} K. K. Kwon, \emph{Some sharp Hodge Laplacian and Steklov eigenvalue estimates for
differential forms}, Calc. Var. Partial Differential Equations
(2016), https://doi.org/10.1007/s00526-016-0977-8.

\bibitem{K} J. R. Kuttler, V. G. Sigillito, \emph{Inequalities for membrane and Stekloff eigenvalues}, J. Math. Anal. Appl. {\bf 23} (1968) 148--160.


\bibitem{LM} S. Li, J. Mao, \emph{Estimates for sums of eigenvalues of the free plate with nonzero Poisson's
ratio}, Proc. Amer. Math. Soc. {\bf 149}(5) (2021)  2167--2177.


\bibitem{LM1} S. Li, J. Mao, \emph{An isoperimetric inequality for a biharmonic Steklov problem}, available online at
arXiv:2105.12401.



\bibitem{RR}  R. Reilly, \emph{On the first eigenvalue of the Laplacian for compact
submanifolds of Euclidean space}, Comment. Math. Helv. {\bf52}(4)
(1977) 525--533.


\bibitem{RS} M. Ritor\'{e}, C. Sinestrari, \emph{Mean Curvature Flow and Isoperimetric
Inequalities}, Adv. Courses Math. CRM Barcelona, Birkh\"{a}user,
Basel, 2010.


\bibitem{sy} R. Schoen, S. T. Yau,  \emph{Lectures on Differential Geometry}, International Press, Cambridge,
2004.

\bibitem{ST} M. W. Stekloff, \emph{Sur les probl\`emes fondamentaux de la physique math\'ematique}, Ann. Sci. \'Ecole Norm. Sup. {\bf 19} (1902) 455--490.

\bibitem{WX1} Q. Wang, C. Xia, \emph{Sharp bounds for the first non-zero Stekloff eigenvalues}, J. Funct. Anal. {\bf 257} (2009) 2635--2654.

\bibitem{WX2} Q. Wang, C. Xia, \emph{Inequalities for the Steklov eigenvalues}, Chaos, Solitons $\&$ Fractals {\bf 48} (2013) 61--67.

\bibitem{WX3} C. Xia, Q. Wang, \emph{Eigenvalues of the Wentzell-Laplace Operator and of the Fourth Order Steklov Problems},
 J. Differential Equations {\bf 264}(10) (2018) 6486--6506.


\bibitem{YWMD} Y. Zhao, C. X. Wu, J. Mao, F. Du, \emph{Eigenvalue comparisons in Steklov eigenvalue problem and some other eigenvalue estimates},
Revista Matem\'{a}tica Complutense {\bf33}(2) (2020) 389--414.

\bibitem{yy} L. Yang, C. Yu, \emph{A higher dimensional generalization of Hersch-Payne-Schiffer
inequality for Steklov eigenvalues}, J.  Funct. Anal. {\bf 272}
 (2017) 4122--4130.


\end{thebibliography}
\end{document}